Solution to Burnside's Problem
S. Bachmuth

1.  **Introduction**

Throughout this paper we fix a prime power $q = p^e$. Unless specifically mentioned otherwise, all groups are 2-generator.

We shall show that the 2-generator Burnside groups of prime power exponent are solvable and thus finite. The main result, which we call the Generalized Burnside Theorem, is a solvability theorem for a more general family of groups. In a final section we outline the extension of these results to k-generator groups.

We begin with a four part theorem. The first part constructs a free group of matrices followed by three of its images. Details and proofs are in sections 2 and 3.

Theorem 1:
(**i**)  We define a (rank 2) free group $F(\mathcal{R}[t, t^{-1}])$ consisting of 2x2 matrices having entries in the Laurent polynomial ring $\mathcal{R}[t, t^{-1}]$, where $\mathcal{R}$ itself is the Laurent polynomial ring with integer coefficients $\mathcal{R} = \mathbb{Z}[x, x^{-1}, y, y^{-1}]$.

(**ii**)  The (homomorphic) image $F(\mathcal{R})$ of $F(\mathcal{R}[t, t^{-1}])$, obtained by setting $t = 1$ in each element of $F(\mathcal{R}[t, t^{-1}])$, is a group of matrices with entries in $\mathcal{R}$ that is isomorphic to the free metabelian group $F/F''$ of rank 2.

(**iii**)  For each positive integer n, there exists a quotient ring $S = S(n)$ of $\mathcal{R}$ such that if n is a prime-power, $F(S)$ is isomorphic to $F/F''F^n$, the free (rank 2) metabelian Burnside group of exponent n.

(**iv**)  $F(S[t, t^{-1}])$ is solvable if (and only if) the integer n is a prime power.

We thus have the following commutative square where the horizontal maps $\alpha$ and $\underline{\alpha}$ come from the ring homomorphism $\mathcal{R} \to S$ and the vertical maps are those determined by sending t to the identity.

$$\begin{array}{ccc} F(\mathcal{R}[t, t^{-1}]) & \xrightarrow{\alpha} & F(S[t, t^{-1}]) \\ \downarrow & & \downarrow \\ F(\mathcal{R}) & \xrightarrow{\underline{\alpha}} & F(S) \end{array}$$

Note that for prime-power exponents, $\underline{\alpha}: F(\mathcal{R}) \to F(S)$ is the Burnside map on the free metabelian group. It is likely also for arbitrary exponents, but our proof is valid only for prime-power exponents.
(The Burnside map on the free group is $\beta: F(\mathcal{R}[t, t^{-1}]) \to G$, where $G \cong F/F^q$.)

Our main result, Theorem GB (Generalized Burnside), requires the concept of *induced map* which is discussed in Section 4. For now, we note that induced map is a familiar concept to many.

Theorem GB : Let $\gamma$ be a mapping of the free group $F(\mathcal{R}[t, t^{-1}])$ onto a group G such that $\gamma$ induces $\underline{\alpha}$. Then G (i.e., the image of $\gamma$) is solvable.

We will define and discuss the relevant facts about induced maps in Section 4. Because the Burnside map on the free group $F(\mathcal{R}[t, t^{-1}])$ induces $\underline{\alpha}$, we have the

Corollary: The Burnside group of exponent q is solvable.

The proof of Theorem GB is divided into two Propositions. Theorem 1(**iv**) is a special case of Proposition 1 where $G = F(\mathcal{S}[t, t^{-1}])$.

Proposition 1: Suppose the surjective map $\gamma : F(\mathcal{R}[t, t^{-1}]) \to G$ sends $\mathcal{R}$ into $\mathcal{S}$. Then G (the image of $\gamma$) is solvable.

Proposition 2: Suppose $\gamma : F(\mathcal{R}[t, t^{-1}]) \to G$ induces $\underline{\alpha}$. Then the entries of $F(\mathcal{R}[t, t^{-1}])$ which are in $\mathcal{R}$ are sent into S.

The Burnside group $F/F^q$ of exponent q has only elements of finite order and is thus, as a consequence of Theorem GB, a finite group. Not all the solvable groups arising from Theorem GB are finite since some possess elements of infinite order. An example of such a group is $F(\mathcal{S}[t, t^{-1}])$. But whether finite or infinite, all of these groups are closely associated with exponent q groups. The assumption that $\gamma$ induces $\underline{\alpha}$ insures that G has elements of order q. A non-obvious implication is that the commutator subgroup of G has exponent q. Indeed, each Generalized Burnside group has a subgroup of exponent q which is strictly larger than the commutator subgroup. Section 5 will illustrate this in the study of $F(\mathcal{S}[t, t^{-1}])$. The presence of such a large subgroup of exponent q in these Generalized Burnside groups lends a most plausible explanation for their solvability.

The proofs of (**i**), (**ii**) and (**iii**) of Theorem 1 are in Section 2,, the proof of part (**iv**) and the more general Proposition 1 is in Section 3. Section 4 contains a discussion of *induced* maps from free groups to free metabelian groups after which the proof of Proposition 2 and hence Theorem GB is given. Section 5 contains a consideration of some properties $F(\mathcal{S}[t, t^{-1}])$. Section 6 considers the extension of Theorem GB to finitely generated groups.

## 2. Metabelian matrix groups

Let F be the free group of rank 2 generated by the 2x2 matrices $M_1$ and $M_2T$ with entries in the Laurent polynomial ring $Z[x, x^{-1}, y, y^{-1}, t, t^{-1}]$.

$$M_1 = \begin{vmatrix} 1 & 1-y \\ 0 & x \end{vmatrix} \qquad M_2T = \begin{vmatrix} yt & 0 \\ 1-xt & 1 \end{vmatrix}$$

That F is free can be seen by taking $x = 1$ and $y = t = -1$ and appealing to a well known theorem. cf. Sanov [9].

Notation: We set $\mathcal{R} = Z[x, x^{-1}, y, y^{-1}]$ and write $F = F(\mathcal{R}[t, t^{-1}])$. The units of $\mathcal{R}$ are the monomials $x^i y^j$ and $-x^i y^j$ where $i, j$ are integers. We call $x^i y^j$ the positive units.

Lemma 1: (i) The image of F upon setting $t = 1$ is isomorphic to the free metabelian group $F/F''$ generated by the matrices $M_1$ and $M_2$ with entries in $\mathcal{R} = Z[x, x^{-1}, y, y^{-1}]$. We put $F(\mathcal{R}) = gp\langle M_1, M_2 \rangle$.

$$M_1 = \begin{vmatrix} 1 & 1-y \\ 0 & x \end{vmatrix} \qquad M_2 = \begin{vmatrix} y & 0 \\ 1-x & 1 \end{vmatrix}$$

(ii) Each element M of $F(\mathcal{R})$ has the form $M = uI + N$, where I is the identity matrix, u is a (positive) unit in $\mathcal{R}$, and the 2x2 matrix N has the form

$$N = \begin{vmatrix} \lambda_1(1-x) & \lambda_1(1-y) \\ \lambda_2(1-x) & \lambda_2(1-y) \end{vmatrix} \quad , \text{ where } \lambda_1, \lambda_2 \text{ in } \mathcal{R} \text{ satisfy}$$

$\lambda_1(1-x) + \lambda_2(1-y) = 1-u$. We will omit the I and write $M = u + N$.

Proof: (i) $F(\mathcal{R})$ is isomorphic to the group of inner automorphisms of the rank 2 free metabelian group as described in [3]. ([3] provides a description of the group of all IA-automorphisms of the free metabelian group which contains the inner automorphisms. For rank 2, all IA-automorphisms of the free metabelian group are inner.)

(ii) A proof is contained in [4], section 2, Proposition(ii). Because of the importance of this Lemma, we shall reproduce the proof in an appendix to this paper.

Notation: Let $v = (1-x, 1-y)$ and write N as above in the form $N = [\lambda_i v]$, where $\lambda_i v$ denotes the $i^{th}$ row of N.

The lower central series for a group G is defined inductively by $G_1 = G$ and for $j > 1$ the $j^{th}$ term is $G_j = [G, G_{j-1}]$.

__Lemma 2:__   If $M = u + N$ is in $F(\mathcal{R})$, where $N = [\lambda_i v]$, then

   (i)   $M^n = u^n + (1 + u + u^2 + \cdots + u^{n-1})N$   for any integer $n > 1$.
   (ii)  If $M$ is in the commutator subgroup of $F(\mathcal{R})$, then $\lambda_1$ and $\lambda_2$ are in $\Sigma$, the augmentation ideal of $\mathcal{R}$. More generally, if $M$ is in the $j^{th}$ term of the lower central series of $F(\mathcal{R})$, then $\lambda_1$ and $\lambda_2$ are in $\Sigma^{j-1}$.

__Proof:__   Note that $N^2 = (1 - u)N$ since $\lambda_1(1-x) + \lambda_2(1-y) = (1 - u)$.

   (i) From Lemma 1(ii) and the fact that $N^2 = (1 - u)N$, we have
$M^2 = (u + N)^2 = u^2 + 2uN + (1-u)N = u^2 + (1+u)N$. Now proceed by induction.
$M^n = M^{n-1}M = \{u^{n-1} + (1 + u + \cdots + u^{n-2})N\}(u + N) = u^n + (1 + u + u^2 + \cdots + u^{n-1})N$ again using $N^2 = (1 - u)N$.

   (ii) Choose a basis for $G = F(\mathcal{R})$ modulo $G_{j+1}$ and, using induction on $j$, compute using this basis. A convenient basis (H. Neumann [8], chapter 3, section 6) is the basic commutators that are not in the second derived group. Notice that if $M$ is in the commutator subgroup, then $u = 1$ and $\lambda_1(1-x) + \lambda_2(1-y) = 0$. Hence $\lambda_1$ and $\lambda_2$ are in the augmentation ideal of $\mathcal{R}$. In fact, $[M_2, M_1] = I + [\lambda_i v]$, where $\lambda_1 = -(1-y)$ and $\lambda_2 = (1-x)$. This is the first case and the start of the induction. It is straightforward to check that if the basic commutator $[M_2, M_1, M_2, ..., M_2, M_1, ..., M_1] = I + [\lambda_i v]$, where $M_1$ and $M_2$ occur a total of $a +1$ and $b +1$ times respectively, then $\lambda_1 = -(1-y)^{b+1}(1-x)^a$ and $\lambda_2 = (1-y)^b(1-x)^{a+1}$. Furthermore, conjugation of a commutator by $M = u + N$ becomes multiplication (of the non identity term) by $u$. Thus, Lemma 2(ii) follows.

Before proceeding we need to introduce more notation.
$q = p^e$ as always denotes a prime power and we let $\phi(q) = p^e - p^{e-1}$, the Euler $\phi$-function. Let $\mathcal{I}(n)$ be the n-cyclotomic ideal, the ideal in $\mathcal{R}$ generated by all elements $1 + u + u^2 + \cdots + u^{n-1}$ where $u$ is a positive unit in $\mathcal{R}$. We denote by $\mathcal{S}$ the quotient ring $\mathcal{R}/\mathcal{I}(n)\Sigma$. In Lemma 3(i)&(iii) we require that n be a prime-power.

__Lemma 3:__ (i)  $\Sigma^{e\phi(q)} \subseteq \mathcal{I}(q)$ in $\mathcal{R}$. Equality holds if (and only if) $e=1$, that is, $q = p^e$ is prime. Furthermore $\Sigma^{e\phi(q)-1} \not\subseteq \mathcal{I}(q)$.
   If $e \geq 2$, $0 \leq j \leq e-1$ and $k \geq e(p^e - p^{e-1}) - j(p^{e-1} - p^{e-2})$, then $p^j\Sigma^k \subseteq \mathcal{I}(q)$.

   (ii) The group generated by the matrices $M_1$ and $M_2$ over $\mathcal{R}/\Sigma^c$ is isomorphic to $F/F''F_c$ the free metabelian nilpotent group of class c.
   (iii) The group $F(\mathcal{S})$ generated by the matrices $M_1$ and $M_2$ over $\mathcal{R}/\mathcal{I}(q)\Sigma$ is isomorphic to $F/F''F^q$, the free metabelian Burnside group of prime power exponent q.

__Proof:__ (i) This is part of Theorem B in [5] and the Proposition on page 241 of [5].

Theorem B is the case $j = 0$ and will suffice for our purposes. However, for computations when $e > 1$, the extra information in [5] is often desirable, even necessary. (Arbitrary finite rank is not addressed here. For now we are concerned only with rank $r = 2$, but we mention that these results from [5] are valid for ranks $r \leq p+1$. For rank $p+2$, Lemma 3 is no longer valid.)

(ii) This follows from Lemmas 1 and 2(ii). For clearly a matrix of $F(\mathcal{R})$ which is in $F(\mathcal{R})_c$ has entries congruent to I modulo $\sum^c$ by Lemma 2(ii). Conversely, suppose $M = u + [\lambda_i v] \equiv I \mod \sum^c$. Then $u = 1$, $\lambda_i v$ are in $\sum^c$ and $\lambda_1(1-x) + \lambda_2(1-y) = 0$. Thus $\lambda_1 = \lambda(1-y)$, $\lambda_2 = -\lambda(1-x)$ and $\lambda$ is in $\sum^{c-2}$. From the proof of Lemma 2(ii), it is clear that M is in $F(\mathcal{R})_c$.

(iii) This follows from Lemmas 1 and 2(i). For again it is clear by Lemma 2(i) that a matrix in $F(\mathcal{R})'$ which is a $q^{th}$ power has entries congruent to I modulo $\mathcal{I}(q)\sum$. Conversely, suppose $M = u + [\lambda_i v] \equiv I$ modulo $\mathcal{I}(q)\sum$. Assume first that M is in the commutator subgroup. Then $u = 1$, $\lambda_i$ are in $\mathcal{I}(q)$, and $\lambda_1(1-x) + \lambda_2(1-y) = 0$. Thus the $\lambda_i$ are in the augmentation ideal and again $\lambda_1 = \lambda(1-y)$, $\lambda_2 = -\lambda(1-x)$ and both $\lambda(1-y)$ and $\lambda(1-x)$ are in $\mathcal{I}(q)$. By Theorem B of [5] or the Theorem of Dark and Newell [6], we conclude that $\lambda(1-y)$ and $\lambda(1-x)$ are in $\sum^{e\phi(q)}$ and thus M is in $F(\mathcal{R})^q$.

Any element M in $F(\mathcal{R})$ can be written $M = M_1^i M_2^j C$, where C is in the commutator subgroup of $F(\mathcal{R})$. Then $M = x^i y^j + [\lambda_i v]$ for suitable $\lambda_i$. By the hypothesis, i and j must each be a multiple of q. Hence, by Lemma 2(i), $M_1^i$ and $M_2^j$ are each congruent to the identity. So we may assume M is in the commutator subgroup. This completes the proof of Lemma 3.

3.   **Power series representations of F ; Proof of Proposition 1**

We return to our free group $F = F(\mathcal{R}[t, t^{-1}]) = gp\langle M_1, M_2T \rangle$ defined at the beginning of the previous section. For computational purposes, we shall describe $F(\mathcal{R}[t, t^{-1}])$ as a power series in $(t-1)^k$, $k = 0,1,2,...$, where the coefficients are 2x2 matrices over $\mathcal{R} = Z[x, x^{-1}, y, y^{-1}]$. Namely, let

$$T = \begin{vmatrix} t & 0 \\ 1-t & 1 \end{vmatrix} \quad \text{and} \quad S = \begin{vmatrix} 1 & 0 \\ -1 & 0 \end{vmatrix}$$

Then,
$$M_2T = M_2(I + (t-1)S) = M_2 + (t-1)M_2S$$
$$(M_2T)^{-1} = (I + (t-1)S)^{-1}M_2^{-1} = (I - (t-1)S + (t-1)^2S^2 - (t-1)^3S^3 + \cdots)M_2^{-1}$$
$$= M_2^{-1} - (t-1)SM_2^{-1} + (t-1)^2SM_2^{-1} - (t-1)^3SM_2^{-1} + \cdots$$

Since $F(\mathcal{R}[t, t^{-1}]) = gp\langle M_1, M_2T \rangle$, we can state formally

Lemma 4: Any element f in $F(\mathcal{R}[t, t^{-1}])$ has a representation of the form
(*) $\qquad f = M_f + (t-1)A_1 + (t-1)^2 A_2 + (t-1)^3 A_3 + \cdots$
where $M_f \in F(\mathcal{R}) = gp\langle M_1, M_2 \rangle$ and the $A_i$ as well as $M_f$ are matrices over $\mathcal{R}$.

We know from Lemma 1(i) that $F/F''$ is the image of F via the map $f \to M_f$. The next few lemmas contain further consequences of this description of $F(\mathcal{R}[t, t^{-1}])$.

Notation: Recall that $\Sigma$ denotes the augmentation ideal of $\mathcal{R}$. If all the entries of a matrix A are in $\Sigma^i$ we say that A is in $\Sigma^i$.

Lemma 5: Suppose $f \in F$ as in (*) is in $F'$ the commutator subgroup of F. Then all the $A_i$ (the coefficients of $(t-1)^i$) are in $\Sigma$.
Proof: In the free group $F(\mathcal{R}[t, t^{-1}])$, set $x = y = 1$ (i.e. factor by the ideal $\Sigma$ in $\mathcal{R}$). What remains is an infinite cyclic group (generated by the matrix we called T at the beginning of this section). Hence, the commutator subgroup of F is in the kernel of this map and the result follows.

Remarks:
(1) This proof readily generalizes when the rank is greater than two. (For rank r, setting $x_1 = \ldots = x_r = 1$ produces a free abelian group of rank r-1. See Section 6.)

(2) If f is in F'', the second derived group of $F(\mathcal{R}[t, t^{-1}])$, then one can show that $A_1$, the coefficient of t-1, is in $\Sigma^3$. However the remaining coefficients are in $\Sigma^2$. For the next lemma, we will be content to merely assume that the coefficients of f in F'' are in $\Sigma$ as asserted in Lemma 5. Using Lemma 5 as a starting point, we show that as we move down the derived series of $F(\mathcal{R}[t, t^{-1}])$, we double the power of the augmentation ideal in the entries of the coefficients of the $(t-1)^i$. The proof of Lemma 6(ii) follows an elegant argument of Vaughan-Lee.

Lemma 6: (i) g is in F'' if and only if g has the form
$$g = I + (t-1)A_1 + (t-1)^2 A_2 + (t-1)^3 A_3 + \cdots$$

(ii) If g is in $F^{(k)}$, the $k^{th}$ derived group of F, $(k > 1)$, and $d = 2^{k-2}$, then
$$g = I + (t-1)^d A_d + (t-1)^{d+1} A_{d+1} + \cdots$$
where $A_d, A_{d+1}, A_{d+2}, \ldots$ are in $\Sigma^d$.
Proof: (i) is a restatement of Lemma 1(i).
  (ii) The proof is by induction on k. The start of our induction is $k = 2$ which is covered by Lemma 5 and part (i) above. Suppose therefore that the result is true for $k \geq 2$ and let $d = 2^{k-2}$. Let g, h be in $F^{(k)}$

$$g = I + (t-1)^d A_d + (t-1)^{d+1} A_{d+1} + \ldots \quad , \quad h = I + (t-1)^d B_d + (t-1)^{d+1} B_{d+1} + \ldots \quad ,$$

and consider $[g, h]$. By induction we assume that $A_i$, $B_i$ are in $\Sigma^d$ for $i \geq d$. Since $g^{-1}, h^{-1}$ are also in $F^{(k)}$, we we can write

$$g^{-1} = I + (t-1)^d C_d + (t-1)^{d+1} C_{d+1} + \ldots \quad , \quad h^{-1} = I + (t-1)^d D_d + (t-1)^{d+1} D_{d+1} + \ldots \quad ,$$

and $C_i$, $D_i$ are in $\Sigma^d$ for $i \geq d$. With the summations over $i \geq d$, let

$$A = \sum (t-1)^i A_i, \quad B = \sum (t-1)^i B_i, \quad C = \sum (t-1)^i C_i, \quad D = \sum (t-1)^i D_i.$$

Then $g\,g^{-1} = (I + A)(I + C) = I + A + C + AC$, and hence $A + C + AC = 0$. Similarly $B + D + BD = 0$. Hence, $[g, h] = (I + A)(I + B)(I + C)(I + D) =$
$I + A + B + C + D + AB + AC + AD + BC + BD + CD + ABC + ABD + ACD + BCD + ABCD =$
$I + AB + AD + BC + CD + ABC + ABD + ACD + BCD + ABCD$.
So if we expand $[g, h]$ as a power series

$$[g, h] = I + (t-1) E_1 + (t-1)^2 E_2 + (t-1)^3 E_3 + \cdots \quad ,$$

then $E_i = 0$ for $i < 2d$, and for $i \geq 2d$, $E_i$ is a linear combination of products of two or more matrices from the set $\{A_d, A_{d+1}, \ldots, B_d, B_{d+1}, \ldots, C_d, C_{d+1}, \ldots, D_d, D_{d+1}, \ldots\}$. By induction, all matrices in this set are in $\Sigma^d$, hence the $E_i$ are in $\Sigma^{2d}$ for $i \geq 2d$. This completes the proof.

Lemma 6(ii) is enough to enable us to show that the group $F(S[t, t^{-1}])$ is solvable. However, this gives a bound for the derived length of $F(S[t, t^{-1}])$ one larger than best possible since we assumed that an element of F'' has all coefficients of $(t-1)^i$ in $\Sigma$, when it is easy to see that they lie in $\Sigma^2$. In order to give the best possible bound, we will use Lemma 7 in place of Lemma 6.

<u>Lemma 7</u>: If $g$ is in $F^{(k)}$, the $k^{th}$ derived group of F, $(k > 1)$, and $d = 2^{k-2}$, then

$$g = I + (t-1)^d A_d + (t-1)^{d+1} A_{d+1} + \cdots$$

where $A_d, A_{d+1}, \ldots$ are in $\Sigma^{2d}$.

Proof: Because of Lemma 6(ii), we only have to show that an element in F'' has all coefficients in $\Sigma^2$. To show this, we use Lemma 5 as a starting point and proceed as in the proof of Lemma 6(ii). The only difference is that the constant term of an element in F' is no longer the identity matrix, but now has the form $u + [\lambda_i v]$, where $\lambda_1$ and $\lambda_2$ are in $\Sigma$. Otherwise the proof proceeds exactly as in Lemma 6(ii). We omit the details.

Recall that the ring $S = S(n) = \mathcal{R}/\mathcal{I}(n)\Sigma$. For the rest of this section we require that $q = p^e$ is a prime power so that we may apply Theorem B of [5] which is needed in Theorem1(**iv**) and Proposition 2. Theorem B links the augmentation and cyclotomic ideals of $\mathcal{R}$. It is summarized in Section 2 (Lemma 3(i)).

<u>Theorem1(**iv**)</u> : If $n = q = p^e$ is a prime power, $F(S[t, t^{-1}])$ is solvable. The derived length of $F(S[t, t^{-1}])$ is at most k where $2^{k-1} \geq e(p^e - p^{e-1}) + 1$.

Proof: We apply Lemma 7 to $F(\mathcal{R}[t, t^{-1}])$ before mapping across to $F(S[t, t^{-1}])$. Choose k large enough so that $2^{k-1} \geq e\phi(q) + 1 = e(p^e - p^{e-1}) + 1$. Then, applying Lemma 3(i), the coefficients of the $(t-1)^i$ are the zero matrix in an expansion of an element of the kth derived group of $F(S[t, t^{-1}])$. Thus the elements in the $k^{th}$ derived group are the identity. This completes the proof of Theorem 1(**iv**).

Comments:
(1) Since $\sum^{e\phi(q)-1} \not\subset \mathcal{I}(q)$, it is easy to see that the value of k in Theorem 1(iv) is best possible. Although $F(S[t, t^{-1}])$ is solvable, a difficult exercise will show that it is not nilpotent. (Factoring by $(t-1)^i$ for suitable i maps $F(S[t, t^{-1}])$ onto a nilpotent group.)
(2) Dark and Newell [6] have shown that the exact nilpotency class of $F/F''F^q$ is $e(p^e - p^{e-1})$ when F has rank 2 and $q = p^e$ is a prime power. Their result can be used to give an alternative method for calculating the solvability class of $F(S[t, t^{-1}])$ when $r = 2$. For our purposes, it is more straightforward to use Theorem B in Bachmuth, Heilbronn, and Mochizuki [5], a method which remains viable for $r \leq p+1$.

We end this section with the first of the two propositions needed for the proof of Theorem GB. Proposition 2 will be considered in the next section.

Proposition 1: Suppose the surjective map $\gamma : F(\mathcal{R}[t, t^{-1}]) \to G$ sends $\mathcal{R}$ into $S$. Then G (the image of $\gamma$) is solvable.

Proof: We already observed that for the special case $G = F(S[t, t^{-1}])$, Proposition 1 is just Theorem 1(**iv**). For this choice of G, the map $F(\mathcal{R}[t, t^{-1}]) \to F(S[t, t^{-1}])$ sends $\mathcal{R}$ to $S$ (while fixing t). For general G, we have as hypothesis that the elements of $\mathcal{R}$ are sent to $S$ in the map $\gamma$. Thus, as in Theorem 1(**iv**), Theorem B in [5] assures that the elements in a large enough power of the augmentation ideal of $\mathcal{R}$ are sent to the zero of S in the map from $F(\mathcal{R}[t, t^{-1}])$ to G. This means that in the expansion of the derived series in G, the elements in a high enough solvability class become the identity in G.

Remarks:
1) In this proof of Proposition 1, we are in effect leaving intact the procedure of Vaughan-Lee as detailed in [1] which showed that $F(S[t, t^{-1}])$ is solvable. We use the hypothesis $\mathcal{R} \to S$ to show that the image of the map $F(\mathcal{R}[t, t^{-1}]) \to G$ is solvable in place of the original map $F(\mathcal{R}[t, t^{-1}]) \to F(S[t, t^{-1}])$, a map defined by $\mathcal{R} \to S$.

2) The image of t determines whether G is finite or infinite. Moreover, it plays a role in the determination of the solvability class. If t is left fixed, then $G = F(S[t, t^{-1}])$ is not only an infinite solvable group, but also has the largest possible solvability class for the given q. Enlarging the ideal will generally decrease the solvability class. In the extreme case of sending t to 1, the ideal becomes as large as possible, and G, which becomes F($S$), has the smallest possible solvability class, (i.e., metabelian or even abelian when q = 2).

3)  Proposition 1 solves a problem in groups (establishing solvability) via calculations in a ring. This procedure was given prominence in Wilhelm Magnus' classic paper [7]. In that paper and in others by various authors since, the terms in a higher commutator series are linked to a 'degree function' in a ring - in our case to the powers of the augmentation ideal. Proposition 1 is achieved by showing that as one proceeds down the derived series of F($R$[t, t$^{-1}$]), the powers of the augmentation ideal of $R$ increase. Theorem B in the joint paper with Heilbronn and Mochizuki [5] ensures that they eventually fall into the q-cyclotomic ideal of $R$ (i.e.,, the zero of $S$ ). Thus, upon the transfer of F($R$[t , t$^{-1}$]) → G, we are able to conclude that G is solvable. This is typical of the type of argument employed by Magnus [7] in his solution of a problem posed by Hopf.

## 4.     Proof of the Generalized Burnside Theorem

<u>Theorem</u> GB: (Generalized Burnside):  Let γ be a mapping of the free group F($R$[t , t$^{-1}$]) onto a group G such that γ  induces  $\underline{\alpha}$ . Then G (i.e., the image of γ)  is solvable.

Since Proposition 1 has been proven, to complete the proof of Theorem GB we only need to establish:

<u>Proposition</u> 2:   Suppose  γ : F($R$[t , t$^{-1}$]) → G  induces  $\underline{\alpha}$ . Then the entries of F($R$[t , t$^{-1}$]) which are in $R$  are sent to S.

Before proving Proposition 2, we review the concept of induced map.  (Section 5 of [1 ] has a slightly more detailed discussion.)

<u>Notation</u>:  For g in F($R$[t , t$^{-1}$]), let g/t  be the matrix obtained by setting t = 1 in g.  Thus g/t  is in F($R$).

<u>Def</u>:  Given  λ: F($R$[t, t$^{-1}$])  →  G.   Define  $\underline{\lambda}$: F($R$)  → G/G″  by  $\underline{\lambda}$(g/t) = λ(g)G″  for g in F($R$[t, t$^{-1}$]).   $\underline{\lambda}$  is the map induced from  λ.

We omit the proof that  $\underline{\lambda}$  is well defined.  A proof is given in Section 5 of [1 ].  The straightforward proof uses the fact that λ takes the second derived group of F($R$[t, t$^{-1}$]) into the second derived group of G.

Remark: It is well known that one can define an induced map on any quotient of a free group by a fully invariant subgroup. For example, automorphisms of free groups and their induced automorphisms (on allowable quotients of free groups) have been much studied. Analogous to homomorphic images, induced maps inherit characteristics of the original. Although information is lost, the concept can prove useful. For induced maps, a critical issue may be the ability to reverse the procedure - i.e., "lift maps".

Our concern here are the maps of $F(\mathcal{R}[t, t^{-1}])$ which induce the Burnside metabelian map $\underline{\alpha}$. Which maps of the free group $F(\mathcal{R}[t, t^{-1}])$ induce $\underline{\alpha}$? The most obvious one is the Burnside map of $F(\mathcal{R}[t, t^{-1}])$. That is, as expected, the Burnside map on the free group induces the Burnside map on the free metabelian group. But there are others, in fact infinitely many maps which induce $\underline{\alpha}$. One of them is $\alpha$ which has $F(S[t, t^{-1}])$ as image. Both $\alpha$ and $\underline{\alpha}$ are defined by sending $\mathcal{R}$ to S. and hence it is not surprising that $\alpha$ induces $\underline{\alpha}$. (A formal proof that $\alpha$ induces $\underline{\alpha}$ can be found in Section 5 of [1].) One can even define $\underline{\alpha}$ as the map induced by $\alpha$. Proposition 2 gives the defining property of these maps. Taken in conjunction with Proposition 1, we conclude that the images of the maps of $F(\mathcal{R}[t, t^{-1}])$ which induce $\underline{\alpha}$ are always solvable groups. Thus, the Generalized Burnside Groups are solvable. In the next section we will learn their relationship to exponent q groups.

Proposition 2. Suppose $\gamma : F(\mathcal{R}[t, t^{-1}]) \to G$ induces $\underline{\alpha}$. Then $\gamma$ sends the entries of $\mathcal{R}$ into S.

Proof: This is Lemma 11 in [1], but we shall give a proof here as well. We begin by noting that since $M_1$ does not have t in any entry, $M_1/t = M_1$. Thus, $M_1$ is in both $F(\mathcal{R})$ and $F(\mathcal{R}[t, t^{-1}]$, and therefore $\gamma(M_1/t) = \gamma(M_1) = \underline{\gamma}(M_1)$. By hypothesis, $\underline{\gamma}(M_1) = \underline{\alpha}(M_1)$, and by definition $\underline{\alpha}(M_1)$ takes the entries of $M_1$ into S. However, the entries of $M_1$ contain the generators of $\mathcal{R}$ and thus $\gamma$ sends all entries of $\mathcal{R}$ into S.

Comments: Proposition 2 by itself imparts much information about G:
(1) The image of the generator $M_1$ of $F(\mathcal{R}[t, t^{-1}])$, being the same as the image of the first generator of the metabelian Burnside map, is thus an element of order q.
(2) Since the elements of $\mathcal{R}$ are mapped to a quotient ring of $\mathcal{R}$, G is a matrix group over a quotient of a Laurent polynomial ring and therefore a Noetherian ring.
(3) Thus far, $\gamma$ behaves exactly like the map of $F(\mathcal{R}[t, t^{-1}])$ to $F(S[t, t^{-1}])$ studied in Section 3. All that remains in defining $\gamma$ is to know what happens to t. But, at this point we can be sure that whatever happens to t, the ideal defining G will contain the ideal defining $F(S[t, t^{-1}])$. Thus G is a homomorphic image of the (solvable) group $F(S[t, t^{-1}])$ and hence G must be solvable. This is therefore a variation of the proof of Proposition 1 given in the previous section.

## 5. The group $G = G(q) = F(S[t, t^{-1}])$.

Having completed the proof of Theorem GB, in this section we investigate some aspects of GB groups.

The group $G = F(S[t, t^{-1}])$ depends on the ring $S = S(n)$ and we can write $G = G(n)$ when we wish to emphasize this dependence on n. As we have seen, when n = q is a prime power, $G = G(q) = F(S[t, t^{-1}])$ is solvable and the results in [5] enable one to determine the exact solvability class. Clearly this class is determined by the ideal in $\mathcal{R}$ that is used to construct $S$, an ideal which depends on q. Many technical questions are easily answered. For example, it is trivial to see that the solvability class of G(q) tends to infinity with increasing q. However, the basic question - why is G(q) solvable for any q? - remains obscure. What is there about G(q) as a group that argues solvability? This question motivates what follows and perhaps shed light on the Burnside Problem.

Assume that $S = S(q)$, where $q = p^e$ is a prime power. The generators of G are the images of $M_1$ and $M_2T$ under $\alpha$. To simplify notation we will use the same letters for the elements of G as for $F(\mathcal{R}[t, t^{-1}])$ but with a line through them. Hence the lined letters are matrices over $S$. Thus

$$G = F(S[t, t^{-1}]) = gp \langle \alpha(M_1), \alpha(M_2T) \rangle = gp \langle \overline{M}_1, \overline{M}_2T \rangle$$

(T is left unchanged since $\alpha$ leaves T unchanged.)

Clearly $\overline{M}_1$ has order q and it is easy to see that $\overline{M}_2T$ has infinite order. Namely, sending x and y to 1 sends $\overline{M}_2T$ to T which is a matrix of infinite order. (Without danger of confusion we are using x and y as elements of $S$ as well as $\mathcal{R}$.) If $\overline{W}$ is an element of G, i.e., a product of $\overline{M}_1$, $\overline{M}_2T$ and their inverses, such that $\overline{W}$ has a nonzero exponential sum in T, then the same proof shows that $\overline{W}$ has infinite order. Next consider the elements with zero exponential sum in T. An example is a conjugate of $\overline{M}_1$ which of course has order q. However, there is a better method to show this, a method which ultimately relies upon Theorem B in [5] (i,e., Lemma 3(i)). That proof shows that any element which has exponential sum zero in T has order a divisor of q. We state this as Theorem 2 and sketch a proof for primes.

Theorem 2: Let $\overline{W}$ be an element of $G = F(S[t, t^{-1}])$.
a) If $\overline{W}$ has exponential sum zero in T, then $\overline{W}^q = 1$.
b) If $\overline{W}$ has a nonzero exponential sum in T, then $\overline{W}$ has infinite order

Proof: Assume q is a prime. If $\overline{W}$ has exponential sum zero in T, then we claim that $\overline{W}^q = 1$. $\overline{W}$ is a product of $\overline{M}_1$, $\overline{M}_2T$ and their inverses. Let W be a preimage of $\overline{W}$ in $F(\mathcal{R}[t, t^{-1}])$. Write the $M_i$ in the form of Lemma 1(ii) of Section 2 and write T in the form $tI + A$, where $tI$ is the scalar matrix with $t$ in the diagonal. As a product of these matrices $W$ has the form $W = uI + V$, where $u = x^i y^j t^k$ and the entries of $V$ are in the augmentation ideal of $\mathcal{R}$. By our assumption, $k = 0$ (the only place we need the hypothesis that the exponential sum of $T$ is zero). Now map $W$ to $\mathcal{W}$, i.e., send $\mathcal{R}$ to $S$. Since $S$ has prime characteristic $q$, $(uI + V)^q = u^q I + V^q$, where $u^q = 1$ and $V^q$ is the zero matrix. Namely, the entries of $V^q$ are in the qth power of the augmentation ideal. By Lemma 3(i), the q-1 power of the augmentation ideal lies in $\mathcal{I}(q)$ and hence the qth power is in $\mathcal{I}(q)\Sigma$ which is the zero of $S$.

We next combine Theorem 1 and Theorem 2. Recall that setting $t = 1$ in a matrix of $G$ sends that matrix into one whose order is a divisor of $q$. (Factoring by the ideal $t - 1$ in $S[t, t^{-1}]$ is factoring by the second derived group in $G$.) Thus as a corollary of Theorems 1 and 2, we have the

Proposition 3: Let $g$ be an element in $G = G(q) = F(S[t, t^{-1}])$. Then either $g^q = 1$ or $g$ has infinite order. If $g$ has infinite order, $g^q \equiv 1$ modulo $G''$.

The ring $S = S(q) = \mathcal{R}/\mathcal{I}(q)\Sigma$ was chosen as precisely the ring which takes the free metabelian group $F(\mathcal{R})$ into $F(S)$, the Burnside metabelian group of exponent $q$. Using the same ring map to change $F(\mathcal{R}[t, t^{-1}])$ into $F(S[t, t^{-1}])$ one might expect to find a substantial number of elements of exponent $q$ in $F(S[t, t^{-1}])$. Theorem 2 describes the ubiquity of these q elements. From this one can hypothesize a convincing (group theoretic) response to the question of *why* $F(S[t, t^{-1}])$ and all the Generalized Burnside groups *are solvable*. Namely, enough elements of exponent q imply (higher) commutator identities as $F(S[t, t^{-1}])$ illustrates. (When q is small, such identities were long known in groups of exponent q.) The ideal in $\mathcal{R}$ used to construct $S$ was chosen to produce elements of exponent q. Then with the help of Theorem B in [5], commutator identities unfold in the Generalized Burnside groups.

Note: As we saw in the previous section, the GB groups are images of $F(S[t, t^{-1}])$.

6. **Arbitrary Generators**

We examine the main ingredients necessary for generalization to k generators. $\mathcal{R} = \mathcal{R}(k) = Z[x_1, x_1^{-1}, \ldots, x_k, x_k^{-1}]$ is the k generator Laurent polynomial ring over the integers. $\mathcal{R}[t, t^{-1}] = Z[x_1, x_1^{-1}, \ldots, x_k, x_k^{-1}, t, t^{-1}]$ is the k+1 generator ring. $\mathcal{I}(q)$ is the q-cyclotomic ideal in $\mathcal{R}$, the ideal generated by all q-cyclotomic elements $1 + u + u^2 + \cdots + u^{q-1}$ where $u$ is a positive unit in $\mathcal{R}$ and q is the prime power $q = p^e$.

$S$ is the quotient ring $\mathcal{R}/\mathcal{I}(q)\Sigma$, where $\Sigma$ is the augmentation ideal of $\mathcal{R}$.

We now describe four k x k matrix groups: $F(\mathcal{R})$, $F(S)$, $F(\mathcal{R}[t, t^{-1}])$ and $F(S[t, t^{-1}])$.

Let v be the k-tuple, $v = (1-x_1, \ldots, 1-x_k)$ and N the k x k matrix $N = [\lambda_i v]$, where $\lambda_i v$ denotes the ith row of N, $i = 1,\ldots,k$. $F(\mathcal{R})$ is the set of all matrices $uI + N$, where I is the identity matrix, u is a (positive) unit in $\mathcal{R}$, and the $\lambda_1, \ldots, \lambda_k$ in $\mathcal{R}$ satisfy $\lambda_1(1-x_1) + \ldots + \lambda_k(1-x_k) = 1-u$. For a matrix M in $F(\mathcal{R})$ we often omit the I and write $M = u + N$. Alternatively, one may define $F(\mathcal{R})$ by the generators:
$$F(\mathcal{R}) = gp\langle M_1, \ldots, M_k \rangle, \text{ where for } j = 1, \ldots, k$$
$M_j = x_j I + [\lambda_i v]$, $\lambda_i = 0$ unless $i = j$ in which case $\lambda_j = 1$.

Replacing $\mathcal{R}$ by $S$ defines $F(S)$ as well as $F(S[t, t^{-1}])$ as soon as we define $F(\mathcal{R}[t, t^{-1}])$. In defining $F(\mathcal{R}[t, t^{-1}])$ we first define k x k matrices $T_i$ ($i = 2,\ldots,k$). $T_i$ has t in the first i-1 diagonal entries, 1 in the remaining diagonal entries, 1-t in the $i^{th}$ row prior to the diagonal, and zeros everywhere else. Notice that $T_2 = T$ when k is 2. It is easy to see that $gp\langle T_2, \ldots, T_k \rangle$ is free abelian of rank k-1, but all that is necessary is that it be abelian (in generalizing Lemma 5 of [1] for arbitrary ranks). We can now define $F(\mathcal{R}[t, t^{-1}])$.
$$F(\mathcal{R}[t, t^{-1}]) = gp\langle M_1, M_2 T_2, M_3 T_3, \ldots, M_k T_k \rangle$$

Those wishing to deal with the Burnside problem for k generators can now follow the earlier procedures for two generator groups. The proof that $F(\mathcal{R})$ is free metabelian of rank k and that $F(S)$ is the free metabelian Burnside group of exponent q and rank k should follow as previously. As before, we have the commutative square where the horizontal maps $\alpha$ and $\underline{\alpha}$ come from the ring homomorphism $\mathcal{R} \to S$ and the vertical maps are those determined by sending t to the identity.

$$\begin{array}{ccc} F(\mathcal{R}[t, t^{-1}]) & \xrightarrow{\alpha} & F(S[t, t^{-1}]) \\ \downarrow & & \downarrow \\ F(\mathcal{R}) & \xrightarrow{\underline{\alpha}} & F(S) \end{array}$$

These preparations should be helpful for those wishing to delve into the k-generator groups. When $k = 2$, our notation reduces to the same as in prior sections.
$\mathcal{R} = \mathcal{R}(2) = Z[x, x^{-1}, y, y^{-1}]$ and $F(\mathcal{R})$, $F(S)$, and $F(\mathcal{R}[t, t^{-1}])$ are as before.
$F(\mathcal{R}) = gp\langle M_1, M_2 \rangle$ and $F(\mathcal{R}[t, t^{-1}]) = gp\langle M_1, M_2 T \rangle$.

# Appendix

<u>Lemma 1(ii)</u>: Each element M of F($\mathcal{R}$) has the form M = uI + N, where I is the identity matrix, u is a positive unit in $\mathcal{R}$, and N = $[\lambda_i v]$ is an r-square matrix whose $i^{th}$ row is $\lambda_i v$. The $\lambda_i$ in $\mathcal{R}$ satisfy $\lambda_1(1-x_1) + \cdots + \lambda_r(1-x_r) = 1 - u$.

<u>Proof</u>: It is a simple matter to check that vM = v for all M in F($\mathcal{R}$). (i.e. verify this for each generator in F($\mathcal{R}$).) Also verify that each generator and it's inverse satisfies Lemma 1(ii). Having verified Lemma 1(ii) for words of length 1, we use induction on the length of a word in F($\mathcal{R}$). Suppose $W_1 W_2$ is an element of F($\mathcal{R}$), where

$W_1 = u_1 + [\lambda_i v]$, $W_2 = u_2 + [\delta_i v]$ and $\lambda_1(1-x_1) + \cdots + \lambda_r(1-x_r) = 1 - u_1$,

$\delta_1(1-x_1) + \cdots + \delta_r(1-x_r) = 1 - u_2$. Then using the fact that vM = v for all M in F($\mathcal{R}$),
$W_1 W_2 = (u_1 + [\lambda_i v])(u_2 + [\delta_i v]) = u_1(u_2 + [\delta_i v]) + [\lambda_i v] = u_1 u_2 + u_1[\delta_i v]) + [\lambda_i v]$.

But $(u_1\delta_1 + \lambda_1)(1-x_1) + \cdots + (u_1\delta_r + \lambda_r)(1-x_r) = u_1(1 - u_2) + (1 - u_1) = 1 - u_1 u_2$.

Thus the result holds for $W_1 W_2$